\documentclass[A4,11pt]{article}
\usepackage{amsfonts}
\usepackage{amsmath,amsthm}
\usepackage{amssymb,amscd}
\usepackage{blkarray}
\usepackage{enumerate}
\usepackage{color}
\usepackage{cite}
\usepackage{esint}
\usepackage{tikz,pgfplots}
\usepackage{pgfplotstable}
\usepackage{tikzscale}
\usepackage[absolute,overlay]{textpos}
\usepackage{multicol}
\usepackage{framed}
\usepackage{lscape}
\usepackage{caption}
\usepackage{subcaption}
\usepackage[inline]{enumitem}
\usepackage{float}

\makeatletter

\newcommand{\Rmnum}[1]{\expandafter\@slowromancap\romannumeral #1@}
\makeatother
\def\Xint#1{\mathchoice
{\XXint\displaystyle\textstyle{#1}}%
{\XXint\textstyle\scriptstyle{#1}}%
{\XXint\scriptstyle\scriptscriptstyle{#1}}%
{\XXint\scriptscriptstyle\scriptscriptstyle{#1}}%
\!\int}
\def\XXint#1#2#3{{\setbox0=\hbox{$#1{#2#3}{\int}$}
\vcenter{\hbox{$#2#3$}}\kern-.5\wd0}}

\def\dashint{\Xint-}

\def\Xint#1{\mathchoice
{\XXint\displaystyle\textstyle{#1}}%
{\XXint\textstyle\scriptstyle{#1}}%
{\XXint\scriptstyle\scriptscriptstyle{#1}}%
{\XXint\scriptscriptstyle\scriptscriptstyle{#1}}%
\!\int}
\def\XXint#1#2#3{{\setbox0=\hbox{$#1{#2#3}{\int}$}
\vcenter{\hbox{$#2#3$}}\kern-.5\wd0}}

\def\dashint{\Xint-}

\theoremstyle{definition}
\newtheorem{theorem}{Theorem}[section]

\newtheorem{proposition}[theorem]{Proposition}

\newtheorem{remark}{Remark}[section]

\begin{document}
\title{Fredholm determinants in the multiparticle hopping asymmetric diffusion model}
\author{\textbf{Eunghyun Lee\footnote{eunghyun.lee@nu.edu.kz}}\\ {\text{Department of Mathematics, Nazarbayev University}}
                                         \date{}   \\ {\text{Kazakhstan}}   }

\date{}
\maketitle

\begin{abstract}
\noindent In this paper we treat  the multiparticle hopping asymmetric diffusion model (MADM) of which initial configuration is such that a single site is occupied by infinitely many particles and all other sites are empty. We show that the probability distribution of the $m^{\textrm{th}}$ leftmost particle's position at time $t$ is represented by a Fredholm determinant.  Also, we consider an exclusion process type model of the MADM, which is the (two-sided) PushASEP. For the  PushASEP with the step Bernoulli initial condition, we find a Fredholm determinant representation of the probability distribution of the $m^{\textrm{th}}$ leftmost particle's position at $t$.
\end{abstract}
\section{Introduction}
\noindent The coordinate Bethe Ansatz has been a quite useful tool to study a certain class of interacting particle systems on the one-dimensional integer lattice $\mathbb{Z}$.  It  provides a direct way to find the transition probabilities of the finite systems. In order to use the coordinate Bethe Ansatz for a certain particle model on $\mathbb{Z}$, the model should meet the \textit{two-particle reducibility}\cite{Korepin}. This requirement can be met by imposing some conditions on the rates of exponential random variables governing the randomness of the model. The asymmetric simple exclusion process (ASEP) is the most extensively studied model. The jumping rates of the particles in the ASEP are constant in each direction of drift. For the ASEP, Tracy and Widom obtained the transition probability of a finite system by using the coordinate Bethe Ansatz and then obtained the probability distribution of the $m^{\textrm{th}}$ leftmost particle's position at time $t$ for an infinite system of which initial configuration is such that all positive integers are occupied and all other sites are empty \cite{TW1}. This probability distribution is represented by a contour integral of a Fredholm determinant \cite{TW2}. After that, they studied the fluctuations of the probability distribution when $t$ tends to $\infty$ for fixed $m$ and when both $m$ and $t$ tend to $\infty$ with fixed ratio of $m$ and $t$ \cite{TW3}. There are also other models that allow the coordinate Bethe Ansatz to be used or that are designed be to \textit{exactly solvable} by the coordinate Bethe Ansatz.  The multiparticle hopping asymmetric diffusion model (MADM) introduced by Sasamoto and Wadati \cite{Sasamoto3} is one of the models exactly solvable by the coordinate Bethe Ansatz (for its generalization, see \cite{Barraquand-Corwin-2016,Povolotsky-2013}).

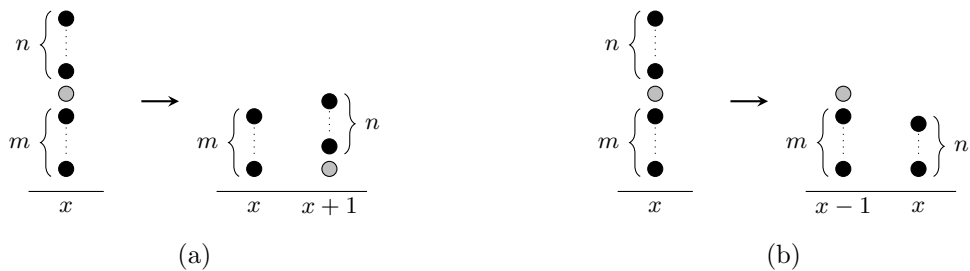
\begin{figure}[h]
\centering
\begin{subfigure}[b]{0.42\textwidth}
\centering
\begin{tikzpicture}[
>=stealth,]
\draw (0,0)--(1,0);
\node at (0.5,-0.2) {\footnotesize $x$};
\draw[fill] (0.5,0.3) circle [radius=0.1];
\draw[dotted] (0.5, 0.5) -- (0.5,0.9);
\draw[fill] (0.5,1) circle [radius=0.1];
\draw[fill=lightgray] (0.5,1.3) circle [radius=0.1];
\draw[fill] (0.5,1.6) circle [radius=0.1];
\draw[dotted] (0.5, 1.8) -- (0.5,2.2);
\draw[fill] (0.5,2.3) circle [radius=0.1];
\draw [decorate,decoration={brace,amplitude= 4 pt}]
(0.3,0.2)--(0.3,1.1) node[midway, left, font=\footnotesize, xshift= -4pt] {$m$};
\draw [decorate,decoration={brace,amplitude=4pt}]
(0.3,1.5)--(0.3,2.4) node[midway, left, font=\footnotesize, xshift=-4pt] {$n$};
\draw [thick, ->] (1.5,1.2)--(2.0,1.2);
\draw (2.5,0)--(4.5,0);
\node at (3,-0.2) {\footnotesize $x$};
\node at (4,-0.2) {\footnotesize $x+1$};
\draw[fill] (3,0.3) circle [radius=0.1];
\draw[dotted] (3, 0.5) -- (3,0.9);
\draw[fill] (3,1) circle [radius=0.1];
\draw[fill=lightgray] (4,0.3) circle [radius=0.1];
\draw[fill] (4,0.6) circle [radius=0.1];
\draw[dotted] (4, 0.8) -- (4,1.1);
\draw[fill] (4,1.2) circle [radius=0.1];
\draw [decorate,decoration={brace,amplitude= 4 pt}]
(2.8,0.2)--(2.8,1.1) node[midway, left, font=\footnotesize, xshift= -4pt] {$m$};
\draw [decorate,decoration={brace,amplitude=4pt,mirror}]
(4.2,0.5)--(4.2,1.3) node[midway, right, font=\footnotesize, xshift=4pt] {$n$};
\end{tikzpicture}
\caption{}
\end{subfigure}
\qquad
\begin{subfigure}[b]{0.42\textwidth}
\centering
\begin{tikzpicture}[
>=stealth,]
\draw (0,0)--(1,0);
\node at (0.5,-0.2) {\footnotesize $x$};
\draw[fill] (0.5,0.3) circle [radius=0.1];
\draw[dotted] (0.5, 0.5) -- (0.5,0.9);
\draw[fill] (0.5,1) circle [radius=0.1];
\draw[fill=lightgray] (0.5,1.3) circle [radius=0.1];
\draw[fill] (0.5,1.6) circle [radius=0.1];
\draw[dotted] (0.5, 1.8) -- (0.5,2.2);
\draw[fill] (0.5,2.3) circle [radius=0.1];
\draw [decorate,decoration={brace,amplitude= 4 pt}]
(0.3,0.2)--(0.3,1.1) node[midway, left, font=\footnotesize, xshift= -4pt] {$m$};
\draw [decorate,decoration={brace,amplitude=4pt}]
(0.3,1.5)--(0.3,2.4) node[midway, left, font=\footnotesize, xshift=-4pt] {$n$};
\draw [thick, ->] (1.5,1.2)--(2.0,1.2);
\draw (2.5,0)--(4.5,0);
\node at (3,-0.2) {\footnotesize $x-1$};
\node at (4,-0.2) {\footnotesize $x$};
\draw[fill] (3,0.3) circle [radius=0.1];
\draw[dotted] (3, 0.5) -- (3,0.9);
\draw[fill] (3,1) circle [radius=0.1];
\draw[fill=lightgray] (3,1.3) circle [radius=0.1];
\draw[fill] (4,0.3) circle [radius=0.1];
\draw[dotted] (4, 0.5) -- (4,0.9);
\draw[fill] (4,0.9) circle [radius=0.1];
\draw [decorate,decoration={brace,amplitude= 4 pt}]
(2.8,0.2)--(2.8,1.1) node[midway, left, font=\footnotesize, xshift= -4pt] {$m$};
\draw [decorate,decoration={brace,amplitude=4pt,mirror}]
(4.2,0.2)--(4.2,1) node[midway, right, font=\footnotesize, xshift=4pt] {$n$};
\end{tikzpicture}
\caption{}
\end{subfigure}
\caption{(a) shows the transition for the case that  the grey particle's random clock with rate $R_n$ rings. (b) shows the transition for the case that  the grey particle's random clock with rate $L_m$ rings. }\label{246-pm-87}
\end{figure}

 The  definition of the MADM is as follows.  Each site $x \in \mathbb{Z}$ is  occupied by $n \in \{0,1,\cdots\}$ particles (or possibly by infinitely many particles) which are stacked upward.   Each particle is equipped with two random clocks exponentially distributed with rates $R_n$ and $L_m$, respectively, where $m$ represents the number of particles below the particle in the stack and $n$ represents the number of particles above the particle in the stack.  If the random clock with rate $R_n$ for a particle at $x$ rings, then  the particle itself and all the $n$ particles in the above move to $x+1$ and form a new stack of particles at $x+1$ (\textit{with particles already occupying $x+1$ if there are}). If the random clock with rate $L_m$ for a particle at $x$ rings, then  the particle itself and all the $m$ particles in the below move to $x-1$ and form a new stack of particles at $x-1$ (\textit{with particles already occupying $x-1$ if there are}).  (See Figure \ref{246-pm-87}.)
\newline
The order of particles in each stack does not matter because all particles are indistinguishable. If the number of particles below or above a given particle changes by particles' jumps, the random clocks reset to have new rates.
 We assume that the rates $R_n$ and $L_m$ are given by
\begin{equation}
 R_n =  p\frac{1-v/u}{1 - \big(v/u\big)^{n+1} }\hspace{0.5cm}\label{RightRate}
  \end{equation}
 and
\begin{equation}
 L_m = q\frac{1-u/v}{1 - \big(u/v\big)^{m+1} } \label{LeftRate}
  \end{equation}
  where $p+q=1=u+v$, $0\leq p \leq 1$ and $0\leq v < u \leq 1$. This condition is required for the integrability  of the model \cite{Ali2,Lee3,Sasamoto3}. A configuration of a finite system with $N$ particles is represented by  $N$-tuples
  \begin{equation*}
  (x_1,\cdots, x_N) \in \mathbb{Z}^N ~~\textrm{with}~x_1 \leq \cdots \leq x_N
  \end{equation*}
  where $x_i$ is the $i^{\textrm{th}}$ leftmost particle's position.   If more than one particle form a stack at a site, we assume the label of the particle in a lower position is less than that of the particle in a higher position. In this paper, we  consider an infinite system of which  configuration can be represented only by
  \begin{equation*}
  (x_1,x_2,\cdots) \in \mathbb{Z}^{\infty} ~~\textrm{with}~x_1 \leq x_2 \leq \cdots.
    \end{equation*}
     Hence, if a configuration of an infinite system is such that $x$ is occupied by infinitely many particles (so, all the sites on the right of $x$ are not occupied and there are only finitely many particles on the left of $x$), then the rate of the random clock for right jump of each particle  at $x$ is equal to $p(1-v/u)$. If this right random clock rings, then the particle and all other infinitely many particles in the above jump to the next right site.
 
If we map $x_i$,  the position of the $i^{\textrm{th}}$ leftmost particle in the MADM, to $x_i +i$, then we obtain an exclusion process type model with \textit{pushing dynamics}, which is called the (two-sided) PushASEP \cite{Ali2,Borodin2,Lee2,Lee3}. (Here, we mean by an exclusion process type model that a site can be occupied by at most one particle in the model.)   The rules of the PushASEP can be stated as follows. Each site $x \in \mathbb{Z}$ is  occupied by at most one particle.  Each particle is equipped with two random clocks exponentially distributed with rates $R_n$ and $L_m$, respectively, where $m$ represents the number of particles contiguously occupying on the left and $n$ represents the number of particles contiguously  occupying on the right.  Here, $R_n$ and $L_m$ are given by (\ref{RightRate}) and (\ref{LeftRate}).   If the random clock with rate $R_n$ for a particle at $x$ rings, then the particle jumps to $x+1$ by pushing all $n$ particles contiguously occupying without gaps on the right by 1 to the right. If the random clock with rate $L_m$ for a particle at $x$ rings, then the particle jumps to $x-1$ by pushing all $m$ particles contiguously occupying without gaps on the left by 1 to the left. (See Figure \ref{309-pm-87}.)
  \begin{figure}[h]
\centering
\begin{subfigure}[b]{0.95\textwidth}
\centering
\begin{tikzpicture}[
>=stealth,scale=0.9]
\draw (-4,0)--(3,0);
\draw[fill] (-3,0.3) circle [radius=0.1];
\draw[fill] (-1.5,0.3) circle [radius=0.1];
\draw[fill] (0.5,0.3) circle [radius=0.1];
\draw[fill] (2,0.3) circle [radius=0.1];
\draw[fill=lightgray] (-0.5,0.3) circle [radius=0.1];
\node at (-3,-0.2) {\footnotesize $x-m$};
\node at (-1.5,-0.2) {\footnotesize $x-1$};
\node at (2,-0.2) {\footnotesize $x+n$};
\node at (0.5,-0.2) {\footnotesize $x+1$};
\node at (-0.5,-0.2) {\footnotesize $x$};
\draw[dotted] (-2.7, 0.3)--(-1.8, 0.3);
\draw[dotted] (1.7, 0.3)--(0.8, 0.3);
\draw [decorate,decoration={brace,amplitude= 4 pt}]
(-3.1,0.5)--(-1.4,0.5) node[midway, above, font=\footnotesize, yshift= 4pt] {$m$};
\draw [decorate,decoration={brace,amplitude= 4 pt}]
(0.4,0.5)--(2.1,0.5) node[midway, above, font=\footnotesize, yshift= 4pt] {$n$};
\draw [thick, ->] (3.5,0.5)--(4,0.5);
\draw (4.5,0)--(12.5,0);
\draw[fill] (5.5,0.3) circle [radius=0.1];
\draw[fill] (7,0.3) circle [radius=0.1];
\draw[fill] (12,0.3) circle [radius=0.1];
\draw[fill] (10,0.3) circle [radius=0.1];
\draw[fill=lightgray] (9,0.3) circle [radius=0.1];
\node at (5.5,-0.2) {\footnotesize $x-m$};
\node at (7,-0.2) {\footnotesize $x-1$};
\node at (12,-0.2) {\footnotesize $x+n+1$};
\node at (10,-0.2) {\footnotesize $x+2$};
\node at (9,-0.2) {\footnotesize $x+1$};
\node at (8,-0.2) {\footnotesize $x$};
\draw[dotted] (5.8, 0.3)--(6.7, 0.3);
\draw[dotted] (10.3, 0.3)--(11.7, 0.3);
\draw [decorate,decoration={brace,amplitude= 4 pt}]
(5.4,0.5)--(7.1,0.5) node[midway, above, font=\footnotesize, yshift= 4pt] {$m$};
\draw [decorate,decoration={brace,amplitude= 4 pt}]
(9.9,0.5)--(12.1,0.5) node[midway, above, font=\footnotesize, yshift= 4pt] {$n$};
\end{tikzpicture}
\caption{}
\end{subfigure}
\\
\begin{subfigure}[b]{0.95\textwidth}
\centering
\begin{tikzpicture}[
>=stealth,scale=0.9]
\draw (-4,0)--(3,0);
\draw[fill] (-3,0.3) circle [radius=0.1];
\draw[fill] (-1.5,0.3) circle [radius=0.1];
\draw[fill] (0.5,0.3) circle [radius=0.1];
\draw[fill] (2,0.3) circle [radius=0.1];
\draw[fill=lightgray] (-0.5,0.3) circle [radius=0.1];
\node at (-3,-0.2) {\footnotesize $x-m$};
\node at (-1.5,-0.2) {\footnotesize $x-1$};
\node at (2,-0.2) {\footnotesize $x+n$};
\node at (0.5,-0.2) {\footnotesize $x+1$};
\node at (-0.5,-0.2) {\footnotesize $x$};
\draw[dotted] (-2.7, 0.3)--(-1.8, 0.3);
\draw[dotted] (1.7, 0.3)--(0.8, 0.3);
\draw [decorate,decoration={brace,amplitude= 4 pt}]
(-3.1,0.5)--(-1.4,0.5) node[midway, above, font=\footnotesize, yshift= 4pt] {$m$};
\draw [decorate,decoration={brace,amplitude= 4 pt}]
(0.4,0.5)--(2.1,0.5) node[midway, above, font=\footnotesize, yshift= 4pt] {$n$};
\draw [thick, ->] (3.5,0.5)--(4,0.5);
\draw (4.5,0)--(12.5,0);
\draw[fill] (5.5,0.3) circle [radius=0.1];
\draw[fill] (7.5,0.3) circle [radius=0.1];
\draw[fill] (12,0.3) circle [radius=0.1];
\draw[fill] (10.5,0.3) circle [radius=0.1];
\draw[fill=lightgray] (8.5,0.3) circle [radius=0.1];
\node at (5.5,-0.2) {\footnotesize $x-m-1$};
\node at (7.5,-0.2) {\footnotesize $x-2$};
\node at (12,-0.2) {\footnotesize $x+n$};
\node at (10.5,-0.2) {\footnotesize $x+1$};
\node at (9.5,-0.2) {\footnotesize $x$};
\node at (8.5,-0.2) {\footnotesize $x-1$};
\draw[dotted] (5.8, 0.3)--(7.2, 0.3);
\draw[dotted] (10.8, 0.3)--(11.7, 0.3);
\draw [decorate,decoration={brace,amplitude= 4 pt}]
(5.4,0.5)--(7.6,0.5) node[midway, above, font=\footnotesize, yshift= 4pt] {$m$};
\draw [decorate,decoration={brace,amplitude= 4 pt}]
(10.4,0.5)--(12.1,0.5) node[midway, above, font=\footnotesize, yshift= 4pt] {$n$};
\end{tikzpicture}
\caption{}
\end{subfigure}
\caption{(a) shows the transition for the case that  the grey particle's random clock with rate $R_n$ rings. (b) shows the transition for the case that  the grey particle's random clock with rate $L_m$ rings. }\label{309-pm-87}
\end{figure}
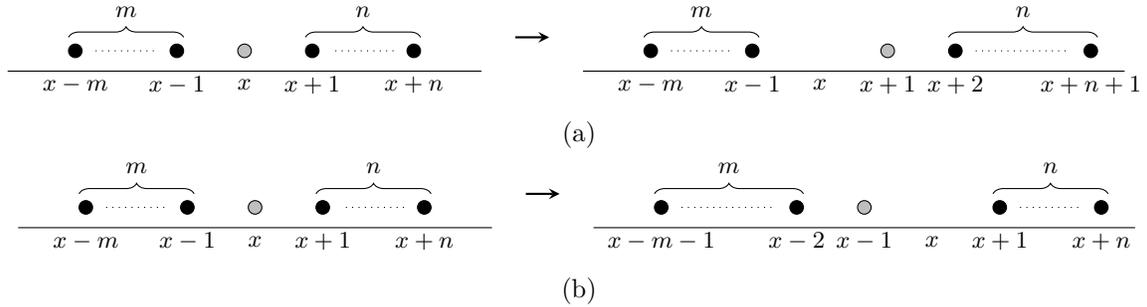
\begin{remark}\textbf{MADM exclusion process and  PushASEP}
The MADM exclusion process is another exclusion process type model which is originated from the MADM (See \cite{Barraquand-Corwin-2016} for its definition). There is a close relationship between the MADM exclusion process and the PushASEP. Indeed, we may think of the holes in the MADM exclusion process as particles in the PushASEP. Then, each transition between configurations in the MADM exclusion process can be mapped to a transition in the PushASEP. But, this mapping is not bijective, that is, there are some transitions in the PushASEP that do not have their corresponding transitions in the MADM exclusion process. For example, let us consider the PushASEP with an initial configuration such that all positive sites are occupied and all other sites are empty (this is called the step initial condition). Here,  each particle in this initial configuration is equipped with a right random clock with rate $p(1-v/u)$ for their right jumps. If the clock for right jump for any particle rings, then the particle jumps to the next right site by pushing  all infinitely many particles on the right by 1. This transition may be considered as a creation of a particle in the MADM with exclusion process but there is no such creation rule in the MADM exclusion process.
\end{remark}

A primary goal of this paper is to find the probability distributions the MADM and  the PushASEP with some special initial conditions. To be more specific, we consider the MADM with an initial configuration such that only one single site is occupied by infinitely many particles and all other sites are empty. For this, we show that the probability distribution of the $m^{\textrm{th}}$ leftmost particle's position at time $t$ is represented by a Fredholm determinant. This will be given in (\ref{1213-am-811}) and (\ref{201-am-811}) in  Section \ref{MADM} . The probability distribution of the same quantity in the PushASEP with the step initial condition is immediately obtained by the mapping explained in the above. We extend this result for the step initial condition to the step Bernoulli initial condition. This will be given in (\ref{417-pm-815}) and (\ref{653-pm-813}) in  Section \ref{Push}. In doing so, we encounter an algebraic identity of which  form is similar to Tracy and Widom's identity in the ASEP with the step Bernoulli initial condition \cite[(9)]{TW4}. Actually, it can be shown that the identity in the PushASEP and the identity in the ASEP are essentially the same. Finally, a comment on asymptotics is given in Remark \ref{1020-am-815}. In the next section, we introduce some notations and previous results for this paper.
\begin{remark}
An initial configuration such that all positive  integers are occupied and all other sites are empty at  $t=0$ in the MADM corresponds to the half flat initial condition in the PushASEP such that all positive even integer sites are occupied and all others are empty. There is an explicit probability formula for the ASEP with the half flat initial condition \cite{Lee} but we report that any trial to find such an explicit formula in a closed form for the PushASEP with the half flat initial condition was unsuccessful.
\end{remark}
\section{Notations and previous results}
As mentioned in Introduction, a configuration of the MADM with $N$ particles is represented by $X=(x_1,\cdots,x_N)\in \mathbb{Z}^N$ with $x_1\leq \cdots \leq x_N$. On the other hand, a configuration of the PushASEP with $N$ particles is $X=(x_1,\cdots,x_N) \in \mathbb{Z}^N$ with $x_1< \cdots < x_N$. When we want to express the position of the $i^{\textrm{th}}$ leftmost particle at time $t$, we write $x_i(t)$. Let $S=\{z_1,\cdots,z_k\} \subset \{1,\cdots, N\}$ with $z_1<\cdots<z_k$, $\sigma(S) = z_1 + \cdots + z_k$, and $|S|$ be the number of elements in $S$. For the parameters $u$ and $v$ in (\ref{RightRate}) and  (\ref{LeftRate}), let
\begin{equation*}
[N] = \frac{u^N - v^N}{u - v},
\end{equation*}
and
\begin{equation*}
[N]! = [N][N-1]\cdots [1], ~{N \brack n} = \frac{[N]!}{[n]![N-n]!}
\end{equation*}
with $[0]!=1.$  The Gaussian binomial coefficient is defined by
\begin{equation*}
 {N \brack n}_{\tau} = \frac{(1-\tau^N)(1-\tau^{N-1})\cdots(1-\tau^{N-n+1})}{(1-\tau)(1-\tau^2)\cdots(1-\tau^{n})},
\end{equation*}
and there is a relationship
\begin{equation*}
{N \brack n} = u^{n(N-n)}{N \brack n}_{\tau}
\end{equation*}
where $\tau = v/u$.
For a complex variable $\xi \neq 0$, let
\begin{equation}\label{energy1}
\varepsilon(\xi) = \frac{p}{\xi} + q\xi -1.
\end{equation}
Let $\mathbb{P}_Y$ be the probability measure of the process that starts  with the initial configuration $Y=(y_1,\cdots, y_N)$. Throughout this paper, we will use the notation $\dashint$ for $\frac{1}{2\pi i}\int$. The author found $\mathbb{P}_Y(x_m(t) = x)$ for the  MADM with an arbitrary deterministic initial configuration \cite[Theorem 3]{Lee3}. Here, we state the formula of $\mathbb{P}_Y(x_m(t) \leq x )$, which is immediately obtained by summing $\mathbb{P}_Y(x_m(t) =x )$ over $x$ from $-\infty$ to $x$.
 \begin{proposition}\cite{Lee3}\label{main2}
Let $0<v < u $,~~ $u+v=1$ and $\tau = v/u$. For the  MADM with $N$ particles, we have
\begin{equation}\label{TaggedProb1}
\begin{aligned}
\mathbb{P}_Y(x_m(t) \leq x )  = &   \sum_{|S| \geq m}(-1)^{m}(uv)^{m(m-1)/2}{|S| - 1 \brack |S|-m}\frac{v^{\sigma(S) - m|S|}}{u^{\sigma(S) - |S|(|S|+1)/2}}\\
 & \times \dashint_{{C}_{R_{z_k}}}\cdots\dashint_{{C}_{R_{z_1}}}
\prod_{ i<j }^k\frac{\xi_{z_i} - \xi_{z_j}}{u + v \xi_{z_i}\xi_{z_j} - \xi_{z_j}} \prod_{i=1}^k\frac{\xi_{z_i}^{x-y_{z_i}}e^{\varepsilon(\xi_{z_i})t}}{ ( 1-\xi_{z_i})}~
d\xi_{z_1}\cdots d\xi_{z_k}
\end{aligned}
\end{equation}
where ${C}_{R_i}$ is a circle centered at 0 with counterclockwise orientation of radius $R_i$ such that $1<R_1 <\cdots< R_N < \tau^{-1}$ and $\varepsilon(\xi_i)$ is given by (\ref{energy1}). The sum runs over all $S\subset \{1,\cdots,N\}$ with $|S|\geq m$.
\end{proposition}
The formula for the same quantity in the PushASEP is as follows.
 \begin{proposition}\cite{Lee3}\label{main2}
Let $0<v < u $,~~ $u+v=1$ and $\tau = v/u$. For the  PushASEP with $N$ particles, we have
\begin{equation}\label{TaggedProb11}
\begin{aligned}
\mathbb{P}_Y(x_m(t) \leq x )  = &   \sum_{|S| \geq m}(-1)^{m}(uv)^{m(m-1)/2}{|S| - 1 \brack |S|-m}\frac{v^{\sigma(S) - m|S|}}{u^{\sigma(S) - |S|(|S|+1)/2}}\\
 & \times \dashint_{{C}_{R_{z_k}}}\cdots\dashint_{{C}_{R_{z_1}}}
\prod_{ i<j }^k\frac{\xi_{z_i} - \xi_{z_j}}{u + v \xi_{z_i}\xi_{z_j} - \xi_{z_j}} \prod_{i=1}^k\frac{\xi_{z_i}^{x-m-(y_{z_i}-z_i)}e^{\varepsilon(\xi_{z_i})t}}{ ( 1-\xi_{z_i})}~
d\xi_{z_1}\cdots d\xi_{z_k}
\end{aligned}
\end{equation}
where $\mathcal{C}_{R_i}$ is a circle centered at 0 with counterclockwise orientation of radius $R_i$ such that $1<R_1 <\cdots< R_N < \tau^{-1}$ and $\varepsilon(\xi_i)$ is given by (\ref{energy1}). The sum runs over all $S\subset \{1,\cdots,N\}$ with $|S|\geq m$.
\end{proposition}
\section{$\mathbb{P}_Y(x_m(t)\leq x)$ with $Y=(0,0,\cdots)$ in the MADM}\label{MADM}
We consider an infinite system of the MADM with the initial configuration $(0,0,\cdots)$. This  can be regarded as the limiting system of the $N$-particle systems with the initial configuration $(0,\cdots,0)$. Hence, we put $y_{z_i} =0$ for each $i$ in (\ref{TaggedProb1}) and then take the limit ($N \to \infty$). A similar approach is also observed in \cite{Lee4}. One may want to symmetrize the integrand, however, to do so,  it should be verified that all (different) contours can be deformed to a single contour for symmetrization.  The following is the algorithm of deforming the contours in (\ref{TaggedProb1}). Let us recall Lemma 2.1 and 2.2 in \cite{Lee3}. We are going to deform $C_{R_1},\cdots,C_{R_{k-1}}$ to $C_{R_k}$, say it $C_{R_k}=C_R$. We can show that the pole of the variable $\xi_{k-1}$ from the $u + v \xi_{i}\xi_{k-1} - \xi_{k-1}~(i<k-1)$ is inside the circle $C_{R_i}$, and the pole of $\xi_{k-1}$ from $u + v \xi_{k-1}\xi_{k} - \xi_{k}$ is outside the contour of $C_{R_k}$.  Hence, we can freely deform the circle  $C_{R_{k-1}}$ to $C_{R_k}$ without encountering any singularities.  Now, we do the same procedure for $C_{R_1},\cdots,C_{R_{k-2}}$ in the sequence of $C_{R_{k-2}}, \cdots, C_{R_1}$. So, we deformed all the contours to $C_{R_k}=C_R$.
\begin{figure}[h]
\centering
\begin{tikzpicture}
\draw (0,0) circle [radius = 1.1];
\draw (0,0) circle [radius =1.7];
\draw (0,0) circle [radius=2];
\node[below] at (0,-1.2) {\footnotesize$C_{R_{k-1}}$};
\node at (1,-2.1) {\footnotesize$C_{R_{k}}$};
\node at (-0.7,0) {\footnotesize$C_{R_{i}}$};
\node at (0.5,0.2) {$\times$};
\node at (2.5,1) {$\times$};
\end{tikzpicture}
\caption{The positions of the poles of $\xi_{k-1}$\label{119-am-811}}
\end{figure}
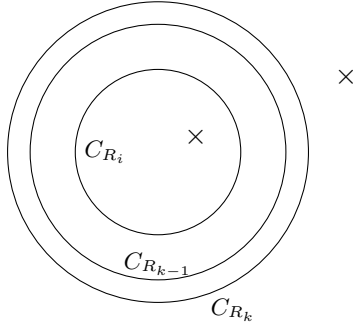
\\
Now, we symmetrize the integrand of  (\ref{TaggedProb1}).  The sum in (\ref{TaggedProb1}) implies $\sum_{k \geq m}\sum_{\substack{S~\textrm{with}~|S| =k}}$.  For a fixed $k = |S|$, the  sum over $S$ with $|S|=k$ in (\ref{TaggedProb1}) is  related to only $v^{\sigma(S)}$ and $u^{\sigma(S)}$. So, the sum over $S$ with $|S|=k$ is written as the sum over all $(z_1,\cdots, z_k)$ with  $1 \leq z_1<z_2<\cdots<z_k\leq N$ because $S=\{z_1,\cdots,z_k\}\subset \{1,\cdots,N\}$. In the limit of  $N \to \infty$, the sum becomes
\begin{equation*}
\sum_{\substack{S~\textrm{with}~|S| =k}}\frac{v^{\sigma(S)}}{u^{\sigma(S)}} = \sum_{1 \leq z_1 < z_2 <\cdots < z_k}\Big(\frac{v}{u}\Big)^{z_1+\cdots+z_k}= \frac{\tau^{k(k+1)/2}}{(1-\tau)(1-\tau^2)\cdots(1-\tau^k)}
\end{equation*}
where $0<\tau = v/u<1$. Let us make variables change as follows: $\xi_{z_i} \to \xi_i$. Now, (\ref{TaggedProb1}) becomes
\begin{equation*}\label{834-pm-89}
\begin{aligned}
\mathbb{P}_Y(x_m(t) \leq x )  = & (-1)^{m}(uv)^{m(m-1)/2}{|S| - 1 \brack |S|-m}\frac{v^{ - m|S|}}{u^{ - |S|(|S|+1)/2}}\frac{\tau^{k(k+1)/2}}{(1-\tau)(1-\tau^2)\cdots(1-\tau^k)} \\
 & \times\int_{{C_R}}\cdots\int_{{C_R}}
\prod_{ i<j }^k\frac{\xi_{i} - \xi_{j}}{u + v \xi_{i}\xi_{j} - \xi_{j}} \prod_{i=1}^k\frac{\xi_{i}^{x}e^{\varepsilon(\xi_{i})t}}{ ( 1-\xi_{i})}~
d\xi_{1}\cdots d\xi_{k} \\
=&  (-1)^{m}(uv)^{m(m-1)/2}{|S| - 1 \brack |S|-m}\frac{v^{ - m|S|}}{u^{ - |S|(|S|+1)/2}}\frac{\tau^{k(k+1)/2}}{(1-\tau)(1-\tau^2)\cdots(1-\tau^k)} \\
 & \times\int_{{C_R}}\cdots\int_{{C_R}}\frac{\prod_{i\neq j}(\xi_j - \xi_i)}{\prod_{i\neq j}(u + v \xi_{i}\xi_{j} - \xi_{j})} \frac{\prod_{i<j }(u + v \xi_{i}\xi_{j} - \xi_{i})}{\prod_{i<j}(\xi_j - \xi_i)}\prod_{i=1}^k\frac{\xi_{i}^{x}e^{\varepsilon(\xi_{i})t}}{ ( 1-\xi_{i})}~
d\xi_{1}\cdots d\xi_{k}.
\end{aligned}
\end{equation*}
Using
\begin{equation}
\sum_{\sigma \in {S}_k} \prod_{ i<j }^k\frac{u + v \xi_{\sigma(i)}\xi_{\sigma(j)} - \xi_{\sigma(i)}}{\xi_{\sigma(j)} - \xi_{\sigma(i)}} = u^{k(k-1)/2}\prod_{i=1}^k\frac{1-\tau^i}{1-\tau} \label{1144-am-810}
\end{equation}
which can be obtained by variable changes  $\xi =\frac{1 +z}{1+z\tau}$ and  (1.4) of Chapter \Rmnum{3} in \cite{MacDonald},  and then  rearranging some terms,  we have
\begin{equation}\label{212-pm-810}
\begin{aligned}
\mathbb{P}_Y(x_m(t) \leq x ) =  &(-1)^{m}\frac{1}{k!}   \sum_{k \geq m}^{\infty} {k - 1 \brack k-m}_{\tau}v^{(k-m)(k-m+1)/2}u^{km+(k-m)(k+m-1)/2}
\\
 & \hspace{0.5cm}\times \int_{{C}_{R}}\cdots\int_{{C}_{R}}
\prod_{ i\neq j }^k\frac{\xi_i - \xi_j}{u + v \xi_i\xi_j - \xi_j}\prod_{i=1}^k\Bigg(\frac{\xi_{i}^{x}e^{\varepsilon(\xi_i)t}}{(1-\tau)(1- \xi_{i})}\Bigg) d\xi_{1}\cdots d\xi_k.
\end{aligned}
\end{equation}
The rest of the procedure of a Fredholm determinant representation for  (\ref{212-pm-810}) is the same as that for the ASEP with the step initial condition. We give some notations. Let
\begin{equation*}
(\lambda;\tau)_m = (1-\lambda)(1-\lambda \tau)\cdots(1-\lambda \tau^{m-1})
\end{equation*}
and  $K$ be an operator on $L^2(C_R)$ defined by
\begin{equation*}
Kf(\xi) = \int_{C_R} K(\xi,\xi')f(\xi')d\xi'
\end{equation*}
where
\begin{equation}\label{1213-am-811}
K(\xi,\xi') = \frac{\xi^{x}e^{\varepsilon(\xi)t}}{u + v \xi\xi' - \xi}\frac{v\xi -u}{1-\tau}.
\end{equation}
Then, just using the same method in  \cite[p.294--295]{TW2}, we can show that (\ref{212-pm-810}) is equal to
\begin{eqnarray}\label{201-am-811}
\mathbb{P}_Y(x_m(t) \leq x ) &=& \int_C\frac{\det(I - \lambda u K)}{(\lambda;\tau)_m}\frac{d\lambda}{\lambda}
\end{eqnarray}
where the contour $C$ is a sufficiently large circle so that all singularities are included in the circle. Comparing (\ref{1213-am-811}) with the kernel of the ASEP in \cite{TW2}, we notice that there is an extra term $\frac{v\xi -u}{1-\tau}$ and $\tau$ in (\ref{201-am-811}) corresponds to $\tau^{-1}$ in \cite{TW2}.
\section{$\mathbb{P}_Y(x_m(t)\leq x)$ with the step Bernoulli initial condition in the PushASEP}\label{Push}
The formula for $\mathbb{P}_Y(x_m(t)\leq x)$ for the PushASEP with initial configuration $(1,2,\cdots)$ is immediately obtained by replacing $x$ by $x-m$ in (\ref{1213-am-811}) from the correspondence between the MADM and the PushASEP.  For an exclusion process type model, the step initial condition is a special case of  the step Bernoulli initial condition \cite{TW4} or the half stationary initial condition \cite{Borodin5}. Here, we  extend $\mathbb{P}_Y(x_m(t)\leq x)$ with the step initial condition in the PushASEP to the case of the step Bernoulli initial condition. We mean the step Bernoulli initial condition by that each positive integer site is occupied with probability $\rho$ with $0<\rho \leq 1$, independently of the other sites at $t=0$, and all other sites are not occupied \cite{TW4}. The derivation below is based on the method in  \cite{TW4}. We first consider the step Bernoulli initial condition on  sites $1,2,\cdots, N$ and then consider its limiting case $N \to \infty$.
\subsection{Averaging the probability measures}
Let $\rho$ be the initial distribution defined by the step Bernoulli initial condition on $1,2,\cdots, N$, and let  $\mathbb{P}_{{\rho},N}$ be the probability measure for the process that starts with $\rho$. Then,
\begin{equation}
\mathbb{P}_{{\rho},N}(x_m(t) \leq x ) = \sum_{Y \subset \{1,\cdots,N\}}\mathbb{P}_Y(x_m(t) \leq x ) \times (\textrm{the probability of $Y$})\label{354-pm-812}
\end{equation}
where $Y$ is the initial distribution with a point mass. It is easy to see that the probability of $Y$ is $\rho^{|Y|}(1-\rho)^{N-|Y|}$. For given two sets $U$ and $V$, let $\sigma(U,V)$ be the number of ordered pairs $(u,v)$ such that $u \in U, v\in V$ and $u \geq v$. Note that $S$ in (\ref{TaggedProb11}) was a subset of $\{1,\cdots,N\}$. After deforming contours, changing variables $\xi_{z_i} \to \xi_i$ and performing some further   manipulations in (\ref{TaggedProb11}), we can show that
\begin{equation}\label{TaggedProb1111}
\begin{aligned}
\mathbb{P}_Y(x_m(t) \leq x )  = &   \sum_{k \geq m}c_k\sum_{\substack{S\subset Y,\\ |S|=k}}\tau^{\sigma(S,Y)} \int_{{C}_{R}}\cdots\int_{{C}_{R}}
\prod_{ i<j }^k\frac{\xi_{i} - \xi_{j}}{u + v \xi_{i}\xi_{j} - \xi_{j}}  \prod_{i=1}^k\frac{\xi_{i}^{x-m}e^{\varepsilon(\xi_{i})t}}{ ( 1-\xi_{i})}\\
 & \times \prod_{i=1}^k{\xi_{i}^{-z_i+\sigma(\{z_i\},Y)}}~
d\xi_{1}\cdots d\xi_{k}
\end{aligned}
\end{equation}
where
\begin{equation*}
c_k = u^{k(k-1)/2}(-1)^{m+1}\tau^{m(m-1)/2}\tau^{-km}{k - 1 \brack k-m}_{\tau}.
\end{equation*}
Here, note that $S$ in  (\ref{TaggedProb1111}) is a subset of $Y$. Substituting (\ref{TaggedProb1111}) into (\ref{354-pm-812}), then
\begin{equation}\label{140-am-812}
\begin{aligned}
& \sum_{k \geq m}c_k\sum_{Y \subset \{1,\cdots,N\}}\sum_{\substack{S\subset Y,\\ |S|=k}}\tau^{\sigma(S,Y)}\rho^{|Y|}(1-\rho)^{N-|Y|} \\
 & \hspace{1cm}\times\int_{{C}_{R}}\cdots\int_{{C}_{R}}
\prod_{ i<j }^k\frac{\xi_{i} - \xi_{j}}{u + v \xi_{i}\xi_{j} - \xi_{j}}  \prod_{i=1}^k\frac{\xi_{i}^{x-m}e^{\varepsilon(\xi_{i})t}}{ ( 1-\xi_{i})} \prod_{i=1}^k{\xi_{i}^{-z_i+\sigma(\{z_i\},Y)}}~
d\xi_{1}\cdots d\xi_{k} \\
=& \sum_{k \geq m}c_k\int_{{C}_{R}}\cdots\int_{{C}_{R}} \sum_{\substack{S \subset \{1,\cdots,N\},\\ |S|=k}}\Bigg(\sum_{\substack{Y\textrm{with}\\ S\subset Y\subset\{1,\cdots,N\}}}\tau^{\sigma(S,Y)}\rho^{|Y|}(1-\rho)^{N-|Y|}\prod_{i=1}^k{\xi_{i}^{\sigma(\{z_i\},Y)}}\Bigg)\\
 & \hspace{1cm}\times\prod_{ i<j }^k\frac{\xi_{i} - \xi_{j}}{u + v \xi_{i}\xi_{j} - \xi_{j}}  \prod_{i=1}^k\frac{\xi_{i}^{x-m}e^{\varepsilon(\xi_{i})t}}{ ( 1-\xi_{i})} \prod_{i=1}^k{\xi_{i}^{-z_i}}~
d\xi_{1}\cdots d\xi_{k}.
\end{aligned}
 \end{equation}
Now, we compute the big parenthesis on right side of (\ref{140-am-812}). As in \cite{TW4}, for $S=\{z_1,\cdots,z_k\}$, let $Y = Y_1\cup Y_2$ where $Y_1 = \{y \in Y : y \leq z_k\}$ and $Y_2 = Y\setminus Y_1$. Then,
\begin{equation}\label{158-am-812}
\begin{aligned} 
&\sum_{\substack{Y\textrm{with}\\ S\subset Y\subset\{1,\cdots, N\}}}\tau^{\sigma(S,Y)}\rho^{|Y|}(1-\rho)^{N-|Y|}\prod_{i=1}^k{\xi_{i}^{\sigma(\{z_i\},Y)}}\\
= & \sum_{\substack{Y\textrm{with}\\ S\subset Y\subset \{1,\cdots,z_k\}}}\tau^{\sigma(S,Y_1)}\rho^{|Y_1|}(1-\rho)^{-|Y_1|}\prod_{i=1}^k\xi_i^{\sigma(\{z_i\},Y_1)} \sum_{Y_2 \subset \{z_{k}+1,\cdots,N\}}\rho^{|Y_2|}(1-\rho)^{N-|Y_2|}
\end{aligned}
\end{equation}
Now, we compute the first sum on the right side in (\ref{158-am-812}) in a similar way to \cite{TW4}. (The second sum on the right side in (\ref{158-am-812}) is $(1-\rho)^{z_k}$ as shown in \cite{TW4}.) If $t_i = z_i -z_{i-1}-1$ and we set $z_0=0$, then $t_i$ implies the number of sites in the gap between $z_i$ and $z_{i-1}$. Let $j_i$ be the number of points of $Y_1$ in the gap between $z_i - z_{i-1}$. Then, the number of $Y_1$ with $j_1,\cdots, j_k$ is ${t_1 \choose j_1}\cdots {t_k \choose j_k}$ and
\begin{equation*}
\sigma(S,Y_1) = \sigma(\{z_1\},Y_1) + \cdots + \sigma(\{z_k\},Y_1)
\end{equation*}
with $\sigma(\{z_i\},Y_1) = j_1+ \cdots+j_i +i$.  Hence, the first sum on the right side of (\ref{158-am-812}) is
\begin{equation*}
\begin{aligned}
&\sum_{j_i \leq t_i}{t_1 \choose j_1}\cdots {t_k \choose j_k}\Big(\frac{\rho}{1-\rho}\Big)^{k + \sum j_i}\tau^{k(k+1)/2 + \sum(k-i+1)j_i}
(\xi_1\cdots\xi_k)^{j_1+1}(\xi_2\cdots\xi_k)^{j_2+1}\cdots \xi_k^{j_k+1}\\
=& \Big(\prod_i\xi_i^i\Big)\tau^{k(k+1)/2}\Big(\frac{\rho}{1-\rho}\Big)^{k}\prod_i\Big(1+\tau^{k-i+1}(\xi_i\cdots\xi_k)\frac{\rho}{1-\rho}\Big)^{t_i}
\end{aligned}
\end{equation*}
by the multi-binomial theorem. Combining this result with $(1-\rho)^{z_k}$ where $z_k$ can be written as $t_1+\cdots +t_k+k$, (\ref{158-am-812}) becomes
\begin{equation}
\Big(\prod_i\xi_i^i\Big)\tau^{k(k+1)/2}\rho^{k}\prod_i\Big(1-\rho+\tau^{k-i+1}(\xi_i\cdots\xi_k){\rho}\Big)^{t_i}.\label{536-pm-812}
\end{equation}
 This is the big parenthesis in (\ref{140-am-812}).
\subsection{Fredholm determinant for the infinite system}
Now, recalling that $z_i = t_1+\cdots +t_i +i$, we see that (\ref{140-am-812}) is
\begin{equation*}
\begin{aligned}
\mathbb{P}_{{\rho},N}(x_m(t) \leq x ) = & \sum_{k \geq m}c_k \sum_{\substack{S \subset \{1,\cdots,N\},\\ |S|=k}}\tau^{k(k+1)/2}\rho^{k} \int_{C_R}\cdots\int_{C_R}\prod_{ i<j }^k\frac{\xi_{i} - \xi_{j}}{u + v \xi_{i}\xi_{j} - \xi_{j}}  \prod_{i=1}^k\frac{\xi_{i}^{x-m}e^{\varepsilon(\xi_{i})t}}{ ( 1-\xi_{i})}\prod_i^k\xi_i^i\\
& \times \prod_{i=1}^k{\xi_{i}^{-z_i+i}}\prod_i^k\Big(1-\rho+\tau^{k-i+1}(\xi_i\cdots\xi_k){\rho}\Big)^{t_i}~
d\xi_{1}\cdots d\xi_{k} \\
 = & \sum_{k \geq m}c_k\tau^{k(k+1)/2}\rho^{k} \int_{C_R}\cdots\int_{C_R}\prod_{ i<j }^k\frac{\xi_{i} - \xi_{j}}{u + v \xi_{i}\xi_{j} - \xi_{j}}  \prod_{i=1}^k\frac{\xi_{i}^{x-m}e^{\varepsilon(\xi_{i})t}}{ ( 1-\xi_{i})}\\
& \times \sum_{\substack{0\leq t_i\textrm{with}\\ \sum t_i \leq N-k}} \prod_i^k\Bigg(\frac{1-\rho+\tau^{k-i+1}(\xi_i\cdots\xi_k){\rho}}{\xi_i\cdots\xi_k}\Bigg)^{t_i}~
d\xi_{1}\cdots d\xi_{k}.
\end{aligned}
\end{equation*}
Since $\xi \in C_R$ with $1<R<\tau^{-1}$ and
\begin{equation*}
\Bigg|\frac{1-\rho+\tau^{k-i+1}(\xi_i\cdots\xi_k){\rho}}{\xi_i\cdots\xi_k}\Bigg|< \big| 1-\rho+\tau^{k-i+1}(\xi_i\cdots\xi_k){\rho}\big| <1,
\end{equation*}
 we have, in the limit $N \to \infty$,
\begin{equation}\label{435-pm-813}
\begin{aligned}
\mathbb{P}_{{\rho}}(x_m(t) \leq x ) :=&\lim_{N \to \infty}\mathbb{P}_{{\rho},N}(x_m(t) \leq x )= \sum_{k \geq m}c_k\tau^{k(k+1)/2}\rho^{k} \int_{C_R}\cdots\int_{C_R}\prod_{ i<j }^k\frac{\xi_{i} - \xi_{j}}{u + v \xi_{i}\xi_{j} - \xi_{j}}  \\
& \times \prod_{i=1}^k\frac{\xi_{i}^{x-m}e^{\varepsilon(\xi_{i})t}}{  1-\xi_{i}}\prod_{i=1}^k\Bigg(\frac{\xi_i\cdots\xi_k}{\xi_i\cdots\xi_k-1+\rho-\tau^{k-i+1}(\xi_i\cdots\xi_k){\rho}}\Bigg)~
d\xi_{1}\cdots d\xi_{k}.
\end{aligned}
\end{equation}
Symmetrizing  the integrand of (\ref{435-pm-813}) is involved with the sum
\begin{equation}\label{500-pm-813}
\begin{aligned}
&\sum_{\sigma \in S_k}\prod_{i<j}\frac{\xi_{\sigma(i)} - \xi_{\sigma(j)}}{u+v\xi_{\sigma(i)}\xi_{\sigma(j)}-\xi_{\sigma(j)}} \frac{\xi_{\sigma(1)}\xi_{\sigma(2)}^2\cdots\xi_{\sigma(k)}^k}
{\prod_{i=1}^k(\xi_{\sigma(i)}\cdots\xi_{\sigma(k)}-1+\rho-\tau^{k-i+1}\rho\xi_{\sigma(i)}\cdots\xi_{\sigma(k)})}
\end{aligned}
\end{equation}
where $\tau = v/u$. This can be simplified by the identity (9) in \cite{TW4} although (\ref{500-pm-813}) looks seemingly different from the quantity in (9) in \cite{TW4}. In (9) in \cite{TW4}, let $p=v$ and $q=u$, and make variable change, that is, $\xi_i \to 1/\xi_i$. Then, we may obtain an equivalent identity
\begin{equation}\label{5000-pm-813}
\begin{aligned}
&\sum_{\sigma \in S_k}\textrm{sgn}(\sigma)\prod_{i<j}\frac{1}{u+v\xi_{\sigma(i)}\xi_{\sigma(j)}-\xi_{\sigma(j)}} \frac{\xi_{\sigma(1)}\xi_{\sigma(2)}^2\cdots\xi_{\sigma(k)}^k}
{\prod_{i=1}^k(\xi_{\sigma(i)}\cdots\xi_{\sigma(k)}-1+\rho-\tau^{k-i+1}\rho\xi_{\sigma(i)}\cdots\xi_{\sigma(k)})}\\
=& u^{k(k-1)/2}\frac{\prod_{i}\xi_i\cdot\prod_{i<j}(\xi_j - \xi_i)}{\prod_i\xi_i - 1 + \rho(1-\tau\xi_i)\cdot\prod_{i\neq j}(u+v\xi_i\xi_j - \xi_j)}.
\end{aligned}
\end{equation}
Hence, by applying (\ref{5000-pm-813}) to (\ref{500-pm-813}), we have
\begin{equation}\label{435-pm-8133}
\begin{aligned}
\mathbb{P}_{{\rho}}(x_m(t) \leq x ) =& \sum_{k \geq m}\frac{1}{k!}u^{k(k-1)/2}c_k\tau^{k(k+1)/2} \int_{C_R}\cdots\int_{C_R}\prod_{ i\neq j }\frac{\xi_{i} - \xi_{j}}{u + v \xi_{i}\xi_{j} - \xi_{j}}  \prod_{i=1}^k\frac{\xi_{i}^{x-m+1}e^{\varepsilon(\xi_{i})t}}{ 1-\xi_{i}}\\
& \times \prod_{i=1}^k\frac{\rho}{\xi_i - 1 + \rho(1-\tau\xi_i)}~
d\xi_{1}\cdots d\xi_{k}.
\end{aligned}
\end{equation}
The rest part of derivation to obtain a Fredholm determinant is trivial because the form of (\ref{435-pm-8133}) is almost the same as the corresponding equation in \cite[p.832]{TW4}. With some observations that our $\tau= v/u$ but $\tau$ in \cite{TW4} is $p/q$, and all coefficients are the same, we finally obtain the Fredholm determinant representation of $\mathbb{P}_{{\rho}}(x_m(t) \leq x )$,
\begin{equation}\label{417-pm-815}
\mathbb{P}_{{\rho}}(x_m(t) \leq x ) = \int_C\frac{\det(I - \lambda K)}{\prod_{j=0}^{m-1}(1-\lambda \tau^j)}\frac{d\lambda}{\lambda}
\end{equation}
where the kernel of $K$ is
\begin{equation}\label{653-pm-813}
K(\xi,\xi') = \frac{\xi^{x-m+1}e^{\varepsilon(\xi)t}}{u+v\xi\xi'-\xi}\cdot\frac{\rho(v\xi-u)}{\xi - 1 + \rho(1-\tau\xi)}
\end{equation} 
and $C$ is sufficiently large so that all the singularities are included. If $\rho=1$, then (\ref{653-pm-813}) reduces to (\ref{1213-am-811}) with $x-m$ instead of $x$, which is the kernel for the PushASEP with the step initial condition. 
\begin{remark}\label{1020-am-815}
In the MADM exclusion process with the step initial condition, it was shown that the first particle's position fluctuates with $t^{1/3}$ scale and the limiting law is governed by the GUE Tracy-Widom distribution \cite{Barraquand-Corwin-2016}. Also, it has been recently shown that the first particle's position in the facilitated exclusion process fluctuates with $t^{1/3}$ scale but its limiting law is governed by GSE Tracy-Widom distribution \cite{Baik-Barraquand-Corwin-Suidan-2017}.   But, in the PushASEP with the step (Bernouli) initial condition, physically, it is expected that the first particle fluctuates with $t^{1/2}$ scale because the first particle will be an asymmetric random walk in the limit $u\to 1$ (but independently of $p$ and $\rho$). (In the limit $u\to 1$, $R_n = p$ for all $n$.) In fact, Tracy and Widom have shown that the first particle's position in the ASEP with the step initial condition fluctuates with $t^{1/2}$ scale and found its liming law, which is non-Gaussian\cite{TW3}. Also,  Lee and Wang have shown that the first particle's position in the $q$-TAZRP with an initial configuration such that the origin is occupied by infinitely many particles and all other sites are empty fluctuates with $t^{1/2}$ and found the same limiting law as Tracy and Widom's for the ASEP \cite{Lee-Wang-2017-arX}. For both the ASEP and the $q$-TAZRP, there is one parameter in the distribution, and when this parameter goes to $0$, the non-Gaussian distribution  ``\textit{converges}" to the Gaussian.  So, it would be interesting if we mathematically confirm $t^{1/2}$ fluctuation and find the limiting law in the PushASEP with the step Bernoulli initial condition where there are three parameters. Also, it would be interesting to see if the limiting distribution is the same as that of the ASEP and the $q$-TAZRP.
\end{remark}
\newpage
\noindent \textbf{Acknowledgement}
\noindent The author would like to thank Guillaume Barraquand and Ivan Corwin for pointing out an error in the earlier version of this paper and providing some valuable comments. This work was partially supported by the Natural Sciences and Engineering Research Council of Canada (NSERC), the Fonds de recherche du Qu\'{e}bec - nature et technologies  (FQRNT), the CRM Laboratoire de Physique Math\'{e}matique, and the social policy grant of Nazarbayev University.


\begin{thebibliography}{99}

\bibitem{Ali2} {Alimohammadi, M., Karimipour, V. and Khorrami, M.:} {A two-parameteric family of asymmetric exclusion processes and its exact solution,} J. Stat. Phys., \textbf{97} 373--394 (1999).
\bibitem{Baik-Barraquand-Corwin-Suidan-2017} {Baik, J., Barraquand, G. and Corwin, I. and Suidan, T.:} Facilitated exclusion process, arXiv:1707.01923.
\bibitem{Barraquand-Corwin-2016} {Barraquand, G. and Corwin, I.:} {The $q$-Hahn asymmetric exclusion process}, Ann. Appl. Probab., \textbf{26}, 2304--2356, (2016).
\bibitem{Borodin5} {Borodin, A., Corwin, I. and Sasamoto, T.:} {From duality to determinants for $q$-TASEP and ASEP},   Ann. Probab., \textbf{42}, 2314--2382, (2014).
\bibitem{Borodin2} {Borodin, A. and Ferrari, P. L.:} {Large time asymptotics of growth models on space-like paths $\Rmnum{1}$:PushASEP}, Electron. J. Probab., \textbf{13} 1380--1418 (2008).
\bibitem{Korepin} {Korepin, V.E., Bogoliubov, and N.M.,Izergin, A.G.:} {Quantum inverse scattering method and correlation functions}, Cambridge University Press (1993).
\bibitem{Lee4} {Korhonen, M., Lee, E.:} {The transition probability and the probability for the left-most particle's position of the $q$-totally asymmetric zero range process}, J. Math. Phys., \textbf{55}, 013301 (2014).
\bibitem{Lee} {Lee, E.:} {Distribution of a particle's position in the ASEP with the alternating initial condition}, J. Stat. Phys., \textbf{140} 635--647 (2010).
\bibitem{Lee2} {Lee, E.:} {Transition probabilites of the Bethe ansatz solvable interacting particle systems,} J. Stat. Phys., \textbf{142} 643--656 (2011).
\bibitem{Lee3} {Lee, E.:} {The current distribution of the multiparticle hopping asymmetric diffusion model,} J. Stat. Phys., \textbf{149} 50--72 (2012).
\bibitem{Lee-Wang-2017-arX} {Lee, E. and Wang, D.:} {Distributions of a particle's position and their asymptotics in the $q$-deformed totally asymmetric zero range process with site dependent jumping rates}, arXiv:1703.08839, (2017).
\bibitem{MacDonald} {MacDonald, I.G.:} {Symmetric functions and Hall polynomials,} Clarendon Press, Oxford (1995).    
\bibitem{Povolotsky-2013} {Povolotsky, A. M.:}  {On the integrability of zero-range chipping models with factorized steady states}, J. Phys. A,  \textbf{46}, 465205, (2013).
\bibitem{Sasamoto3} {Sasamoto, T. and Wadati, M.:} {One-dimensional asymmetric diffusion model without exclusion}, Phys. Rev. E, \textbf{58} 4181--4190 (1998).
\bibitem{TW1} {Tracy, C. A. and Widom, H.:} {Integral formulas for the asymmetric simple exclusion process,} Commun. Math. Phys., \textbf{279} 815--844 (2008), {Erratum : Commum. Math. Phys. \textbf{304}, 875--878 (2011)}.
\bibitem{TW2} { Tracy, C. A. and  Widom, H.:} {A Fredholm determinant representation in ASEP,} {J. Stat. Phys.}, \textbf{132} 291--300 (2008).
\bibitem{TW3} { Tracy, C. A. and  Widom, H.:} {Asymptotics in ASEP with step initial condition,} Commun. Math. Phys., \textbf{290} 129--154 (2009).
\bibitem{TW4} { Tracy, C. A. and  Widom, H.:} {On ASEP with step Bernoulli initial condition}, J. Stat. Phys., \textbf{137} 825--838 (2009).


\end{thebibliography}
\end{document}